\documentclass[12pt]{amsart}

\usepackage{amssymb}
\usepackage{graphicx}
\usepackage{xypic}
\usepackage{hyperref}

\newtheorem{theorem}{Theorem}[section]
\newtheorem{definition}[theorem]{Definition}

\newtheorem{proposition}[theorem]{Proposition}
\newtheorem{corollary}[theorem]{Corollary}

\newtheorem{lemma}[theorem]{Lemma}

\newcommand{\Ker}{\mathop{\mathrm{Ker}}}
\newcommand{\End}{\mathop{\mathrm{End}}}
\newcommand{\Hom}{\mathop{\mathrm{Hom}}}
\newcommand{\Span}{\mathop{\mathrm{span}}}
\newcommand{\diag}{\mathop{\mathrm{diag}}}
\newcommand{\C}{\mathbb C}
\newcommand{\Z}{\mathbb Z}
\newcommand{\N}{\mathbb N}
\newcommand{\R}{\mathbb R}
\newcommand{\T}{\mathbb T}

\begin{document}

\title{Invariants of the half-liberated orthogonal group}

\author{Teodor Banica} \address{T.B.: Department of Mathematics, Toulouse 3
  University, 118 route de Narbonne, 31062 Toulouse, France. {\tt
    banica@math.ups-tlse.fr}}

\author{Roland Vergnioux} \address{R.V.: Department of Mathematics, Caen
  University, BP 5186, 14032 Caen Cedex, France. {\tt
    vergnioux@math.unicaen.fr}}


\begin{abstract}
  The half-liberated orthogonal group $O_n^*$ appears as intermediate quantum
  group between the orthogonal group $O_n$, and its free version $O_n^+$. We
  discuss here its basic algebraic properties, and we classify its irreducible
  representations. The classification of representations is done by using a
  certain twisting-type relation between $O_n^*$ and $U_n$, a non abelian
  discrete group playing the role of weight lattice for $O_n^*$, and a number of
  methods inspired from the theory of Lie algebras. We use these results for
  showing that the discrete quantum group dual to $O_n^*$ has polynomial growth.
\end{abstract}

\maketitle

\section*{Introduction}

The quantum groups introduced by Drinfeld in \cite{dri} have played a prominent
role in various areas of mathematics and physics. In addition to Drinfeld's
discovery, Woronowicz's axiomatization in \cite{wo2}, \cite{wo3} of the compact
quantum groups has been very influential as well and opened the way for the
search of new examples. In particular it allowed the discovery by Wang of the
free quantum groups \cite{wa1}, which have been subject of several systematic
investigations. Since then other families of examples have been discovered,
building up a fast evolving area:
\begin{enumerate}
\item The first two quantum groups are $O_n^+$, $U_n^+$, introduced in
  \cite{wa1}. These led to a number of general developments, including the study
  of connections with subfactors, noncommutative geometry and free probability
  \cite{ba1}, \cite{bco}, \cite{gos} and a number of advances in relation with
  operator algebras \cite{pro}, \cite{vve}, \cite{ve1}.
\item The third quantum group is $S_n^+$, introduced in \cite{wa2}. This led to
  the quite amazing world of quantum permutation groups, heavily investigated in
  the last few years. These quantum groups allowed in particular a clarification
  of the relation with subfactors \cite{bbi}, noncommutative geometry \cite{bgs}
  and free probability \cite{ksp}.
\item The fourth quantum group is $H_n^+$, recently introduced in
  \cite{bbc}. This quantum group gave rise as well to a number of new
  investigations, which are currently under development. Let us mention here the
  opening world of quantum reflection groups \cite{bv2}, and the new formalism
  of easy quantum groups \cite{bsp}.
\end{enumerate}

\vspace{3ex}\pagebreak

In this paper we study the fifth ``new'' quantum group, namely the
half-liberated orthogonal group $O_n^*$, constructed in \cite{bsp}. We believe
that, as it was the case with its predecessors $O_n^+$, $U_n^+$, $S_n^+$,
$H_n^+$, this quantum group will open up as well a new area: namely, that of the
``root systems'' for the compact quantum groups.

The quantum group $O_n^*$ is constructed as follows. Consider the basic
coordinate functions $u_{ij}:O_n\to\C$. These commute with each other, form an
orthogonal matrix, and generate the algebra $C(O_n)$. Wang's algebra $C(O_n^+)$
is obtained by simply removing the ``commutativity'' assumption. As for
obtaining the ``half-liberated'' algebra $C(O_n^*)$, the commutativity condition
$ab=ba$ with $a$, $b\in\{u_{ij}\}$ should be replaced by the weaker condition
$abc=cba$, for $a$, $b$, $c\in\{u_{ij}\}$.

Summarizing, the quantum group $O_n^*$ appears via a kind of tricky
```weakening'' of Wang's original relations in \cite{wa1}. Observe that we have
$O_n\subset O_n^*\subset O_n^+$.

One can prove that, under a suitable ``easiness'' assumption, $O_n^*$ is the
only quantum group between $O_n$ and $O_n^+$. This abstract result, to be proved
in this paper, justifies the name ``half-liberated'', and provides a first
motivation for the study of $O_n^*$. In fact $O_n^*$ was introduced in as one of
the 15 easy intermediate quantum groups $S_n\subset G\subset O_n^+$.

\vspace{3ex}

In this paper we perform a systematic study of $O_n^*$, and in particular of its
category of representations. After discussing the first features of the
definition, we describe in Sections~2~--~4 some $\Hom$-spaces of this category
in terms of Brauer diagrams and derive two consequences: a connection with the
group $U_n$ via the projective version $PO^*_n$, and the uniqueness result
mentionned above.

Then we undertake the classification of irreducible representations of
$O_n^*$. The main technical novelty is the use of diagonal groups and root
systems in a quantum framework. Diagonal groups are meant to be replacements for
maximal torii in good situations and are introduced in Section~5. In the case of
$O_n^*$ the diagonal group provides a noncommutative weight lattice which is
used in Sections~6 and 7 together with a subtle relation to the classical group
$U_n$ to classify representations of $O_n^*$.

Finally we derive in Sections~8, 9 some applications of the classification of
irreducible representations to fusion rules, Cayley graph and growth. Let us
mention here a quite surprising feature of our results: although $O_n^*$ is an
intermediate subgroup between two orthogonal groups $O_n$, $O_n^+$ with
commutative fusion rules, its fusion rules are noncommutative and its exponent
of polynomial growth is the same as for $SU_n$. This shows also that $O_n^*$ is
not monoidally equivalent, in the sense of \cite{brv}, to any known compact
quantum group so far, in particular it is the first original example of a
compact quantum group with exponential growth as considered in \cite{bv1}.

\section{Half-liberation}

Given a compact group $G$, the algebra of complex continuous functions $C(G)$ is
a Hopf algebra, with comultiplication, counit and antipode given by:
\begin{eqnarray*}
  \Delta(\varphi)&=&((g,h)\to\varphi(gh))\\
  \varepsilon(\varphi)&=&\varphi(1)\\
  S(\varphi)&=&(g\to\varphi(g^{-1}))
\end{eqnarray*}

Consider in particular the orthogonal group $O_n$. This is a real algebraic
group, and we denote by $x_{ij}:O_n\to\R$ its basic coordinates,
$x_{ij}(g)=g_{ij}$.

The matrix $x=(x_{ij})$ is by definition orthogonal, in the sense that all its
entries are self-adjoint, and we have $xx^t=x^tx=1$. Moreover, it follows from
the Stone-Weierstrass theorem that the elements $x_{ij}$ generate $C(O_n)$ as a
$C^*$-algebra.

These observations lead to the following presentation result.

\begin{theorem}
  $C(O_n)$ is the universal unital commutative $C^*$-algebra generated by the
  entries of an $n\times n$ orthogonal matrix $x$. The maps given by
  \begin{eqnarray*}
    \Delta(x_{ij})&=&\sum x_{ik}\otimes x_{kj}\\
    \varepsilon(x_{ij})&=&\delta_{ij}\\
    S(x_{ij})&=&x_{ji}
  \end{eqnarray*}
  are the comultiplication, counit and antipode of $C(O_n)$.
\end{theorem}

\begin{proof}
  The first assertion is a direct application of the classical theorems of
  Stone-Weierstrass and Gelfand. The second assertion follows by transposing the
  usual rules for the matrix multiplication, unit and inversion.
\end{proof}

The following key definition is due to Wang \cite{wa1}.

\begin{definition}
  $A_o(n)$ is the universal unital $C^*$-algebra generated by the entries of an
  $n\times n$ orthogonal matrix $u$. The maps given by
  \begin{eqnarray*}
    \Delta(u_{ij})&=&\sum u_{ik}\otimes u_{kj}\\
    \varepsilon(u_{ij})&=&\delta_{ij}\\
    S(u_{ij})&=&u_{ji}
  \end{eqnarray*}
  are the comultiplication, counit and antipode of $A_o(n)$.
\end{definition}

It is routine to check that $A_o(n)$ satisfies the general axioms of Woronowicz
in \cite{wo2}. This tells us that we have the heuristic formula
$A_o(n)=C(O_n^+)$, where $O_n^+$ is a certain compact quantum group, called free
version of $O_n$. See \cite{wa1}.

It is known that we have $A_o(2)=C(SU_2^{-1})$. More generally, the algebra
$A_o(n)$ with $n\geq 2$ arbitrary shares many properties with the algebra
$C(SU_2)$. See \cite{ba1}, \cite{bco}.

The following definition is from the recent paper \cite{bsp}.

\begin{definition}
  The half-liberated orthogonal quantum algebra is
$$A_o^*(n)=A_o(n)/\left\langle abc=cba\,\big|\,a, b, c\in \{u_{ij}\}\right\rangle $$
with comultiplication, counit and antipode coming from those of $A_o(n)$.
\end{definition}

It is routine to check that the comultiplication, counit and antipode of
$A_o(n)$ factorize indeed, and that $A_o^*(n)$ satisfies the general axioms of
Woronowicz in \cite{wo2}. 

In order to get some insight into the structure of $A_o^*(n)$, we first examine
its ``cocommutative version''. We have the following analogue of Definition~1.3.

\begin{definition}
  We consider the discrete group
  $$L_n=\Z_2^{*n}/\left\langle abc=cba\,\big|\, a, b, c\in \{g_i\}\right\rangle$$
  where $g_1, \ldots, g_n$ with $g_i^2=1$ are the standard generators of
  $\Z_2^{*n}$.
\end{definition}

Observe that we have surjective group morphisms $\Z_2^{*n}\to
L_n\to\Z_2^n$. As shown in Proposition~1.6 below, these morphisms are not
isomorphisms in general.

We recall that for any discrete group $\Gamma$, the group algebra $C^*(\Gamma)$
is a Hopf algebra, with comultiplication, counit and antipode given by:
\begin{eqnarray*}
  \Delta(g)&=&g\otimes g\\
  \varepsilon(g)&=&1\\
  S(g)&=&g^{-1}
\end{eqnarray*}

The interest in the above group $L_n$ comes from the following result.

\begin{proposition}
  We have quotient maps as follows:
  $$\xymatrix{A_o(n) \ar[d] \ar[r] & A_o^*(n) \ar[d] \ar[r] & C(O_n) \ar[d] \\
  C^*(\Z_2^{*n}) \ar[r] & C^*(L_n) \ar[r] & C^*(\Z_2^n)}$$
\end{proposition}

\begin{proof}
  The vertical maps can be defined indeed by $u_{ij}\to\delta_{ij}g_i$, by using
  the universal property of the algebras on top. Observe that these maps are
  indeed Hopf algebra morphisms, because the formulae of $\Delta$, $\varepsilon$, $S$
  in Definition~1.2 reduce to the above cocommutative formulae, after performing
  the identification $u_{ij}=0$ for $i\neq j$.
\end{proof}

The group $L_n$ will play an important role in the present paper in relation
with the representation theory of $O^*_n$. In particular we will obtain in
section 6 an abstract isomorphism $L_n \simeq \Z^{n-1} \rtimes \Z_2$. Let us
start with a simpler statement that we use for Theorem~1.7. Note that for $n=2$
this already gives $L_2 = D_\infty = \Z \rtimes \Z_2$.

\begin{proposition}
  The groups $L_n$ are as follows:
  \begin{enumerate}
  \item At $n=2$ we have $\Z_2^{*2}=L_2\neq\Z_2^2$.

  \item At $n\geq 3$, we have $\Z_2^{*n}\neq L_n\neq\Z_2^n$.
  \end{enumerate}
\end{proposition}

\begin{proof}
  We denote by $g$, $h$ the standard generators of $\Z_2^{*2}$.

  (1) We know that $L_2$ appears as quotient of $\Z_2^{*2}$ by the
  relations $abc=cba$, with $a$, $b$, $c\in\{g, h\}$. In the case $a=b$ or $b=c$ this
  is a trivial relation (of type $k=k$), and in the case $a\neq b$, $b\neq c$ we
  must have $a=c$, so once again our relation is trivial (of type
  $klk=klk$). Thus we have $L_2=\Z_2^{*2}$, which gives the result.

  (2) The first assertion is clear, because the equality $g_1g_2g_3=g_3g_2g_1$
  doesn't hold in $\Z_2^{*n}$. Observe now that we have a quotient map
  $L_3\to\Z_2^{*2}$, given by $g_1\to g$, $g_2\to g$, $g_3\to h$. This
  shows that $L_3$ is not abelian. We deduce that $L_n$ with $n\geq 3$ is not
  abelian either, so in particular it is not isomorphic to $\Z_2^n$.
\end{proof}

\begin{theorem}
  The algebras $A_o^*(n)$ are as follows:
  \begin{enumerate}
  \item At $n=2$ we have $A_o(2)=A_o^*(2)\neq C(O_2)$.

  \item At $n\geq 3$, we have $A_o(n)\neq A_o^*(n)\neq C(O_n)$.
  \end{enumerate}
\end{theorem}

\begin{proof}
  The three non-equalities in the statement follow from the three non-equalities
  in Proposition~1.6. In order to prove the remaining statement
  $A_o^*(2)=A_o(2)$, consider the fundamental corepresentation of $A_o(2)$:
$$u=\begin{pmatrix}x&y\\ z&t\end{pmatrix}$$

The elements $x$, $y$, $z$, $t$ are by definition self-adjoint, and satisfy the relations
making $u$ unitary. These unitary relations can be written as follows:
\begin{eqnarray*}
  x^2&=&t^2\\
  y^2&=&z^2\\
  x^2+y^2&=&1\\
  xy+zt&=&0\\
  xz+yt&=&0
\end{eqnarray*}

With these relations in hand, the verification of the relations of type
$abc=cba$ with $a$, $b$, $c\in\{x, y, z, t\}$ is routine, and we get
$A_o^*(2)=A_o(2)$.
\end{proof}

\section{Brauer diagrams}

In this section and in the next one we discuss some basic properties of
$A_o^*(n)$, by ``interpolating'' between some well-known results regarding
$A_o(n)$ and $C(O_n)$.

For $k$, $l$ with $k+l$ even we consider the pairings between an upper sequence
of $k$ points, and a lower sequence of $l$ points. We make the convention that
the $k+l$ points of a pairing $p$ are counted counterclockwise, starting from
bottom left. As an example we draw hereafter the diagram corresponding to the
pairing $\{\{1,9\},\{2,4\},\{3,7\},\{5,8\},\{6,10\}\}$ for $k=l=5$. These
pairings, also called Brauer diagrams, are taken as usual up to planar
isotopy. See e.g. \cite{wen}.

\setlength{\unitlength}{1pt}
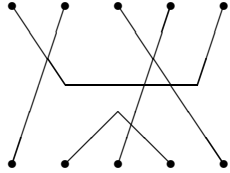
\begin{figure}[h!]
  \centering
  \begin{picture}(100,100)
    \put(0,60){\circle*{3}}\put(20,60){\circle*{3}}
    \put(40,60){\circle*{3}}\put(60,60){\circle*{3}}
    \put(80,60){\circle*{3}}\put(80,0){\circle*{3}}
    \put(60,0){\circle*{3}}\put(40,0){\circle*{3}}
    \put(20,0){\circle*{3}}\put(0,0){\circle*{3}}
    \put(20,0){\line(1,1){20}}
    \put(60,0){\line(-1,1){20}}
    \put(40,0){\line(1,3){20}}
    \put(40,60){\line(2,-3){40}}
    \put(20,60){\line(-1,-3){20}}
    \put(0,60){\line(2,-3){20}}
    \put(80,60){\line(-1,-3){10}}
    \put(20,30){\line(1,0){50}}
  \end{picture}
  \caption{An element of $P(5,5)$ not in $E(5,5)$.}
\end{figure}

\begin{definition}
  We use the following sets of partitions:
  \begin{enumerate}
  \item $P(k, l)$: all pairings.
  \item $E(k, l)$: all pairings with each string having an even number of
    crossings.
  \item $N(k, l)$: all pairings having no crossing at all.
  \end{enumerate}
\end{definition}

The partitions in $N(k, l)$ are familiar objects, also called Temperley-Lieb
diagrams. Observe that the number of crossings for each string of a pairing is
invariant under planar isotopy, so the middle set $E(k, l)$ is indeed
well-defined.

We make the convention that for $k+l$ odd the above three sets are defined as
well, as being equal to $\emptyset$. Observe that for any $k$, $l$ we have
embeddings as follows:
$$N(k, l)\subset E(k, l)\subset P(k, l)$$

\begin{proposition}
  For $p\in P(k, l)$, the following are equivalent:
  \begin{enumerate}
  \item Each string has an even number of crossings (i.e. $p\in E(k, l)$).
  \item The number of points between the two legs of any string is even.
  \item When labelling the points $ababab...$, each string joins an ``$a$'' to a
    ``$b$''.
  \end{enumerate}
\end{proposition}

\begin{proof}
  The equivalence between (2) and (3) is clear. Now fix a string $s$ of $p$ and
  make the following remark: a string $t \neq s$ crosses $s$ {\bf iff} it has
  exactly $1$ end between the two legs of $s$ ; otherwise it has $0$ or $2$
  ends between these legs. This shows the equivalence between (1) and (2).
\end{proof}

The interest in the Brauer diagrams comes from the fact that they encode several
key classes of linear maps, according to the following construction.

\begin{definition}
  Associated to any partition $p\in P(k, l)$ and any $n\in\N$ is the
  linear map $T_p:(\C^n)^{\otimes k}\to (\C^n)^{\otimes l}$ given
  by
  $$T_p(e_{i_1}\otimes\ldots\otimes e_{i_k})=\sum_{j_1\ldots j_l}\delta\begin{pmatrix}i_1
    &\ldots&i_k\\ &p&\\ j_1&\ldots&j_l\end{pmatrix} e_{j_1}\otimes\ldots\otimes
  e_{j_l}$$ 
  where $e_1, \ldots, e_n$ is the standard basis of $\C^n$, and the
  $\delta$ symbol is defined as follows: $\delta=1$ if each string of $p$ joins
  a pair of equal indices, and $\delta=0$ if not.
\end{definition}

Here are a few examples of such linear maps, which are of certain interest for
the considerations to follow:
\begin{eqnarray*}
  T_{\vert\vert}(e_a\otimes e_b)
  &=&e_a\otimes e_b\\
  T_{\slash\hskip-1.5mm\backslash}(e_a\otimes e_b)
  &=&e_b\otimes e_a\\
  T_{\overset{\scriptscriptstyle\sqcup}{\scriptscriptstyle\sqcap}}(e_a\otimes e_b)
  &=&\delta_{ab}\textstyle\sum_ce_c\otimes e_c
\end{eqnarray*}

It is known since Brauer that the linear maps $T_p$ with $p$ ranging over all
diagrams in $P(k, l)$ span the tensor category associated to $O_n$. See
\cite{bra}, \cite{csn}.

\begin{theorem}
  We have the following results:
  \begin{enumerate}
  \item For $C(O_n)$ we have $\Hom(u^{\otimes k}, u^{\otimes l})=\Span(T_p~|~p\in
    P(k, l))$.
  \item For $A_o^*(n)$ we have $\Hom(u^{\otimes k}, u^{\otimes l})=\Span(T_p~|~p\in
    E(k, l))$.
  \item For $A_o(n)$ we have $\Hom(u^{\otimes k}, u^{\otimes l})=\Span(T_p~|~p\in
    N(k, l))$.
  \end{enumerate}
\end{theorem}

\begin{proof}
  The first assertion is Brauer's theorem, see Section~5.a) of \cite{bra}. The
  third assertion is proved in \cite{ba1}, see Proposition~2 and the following
  Remarque there. The middle assertion is Theorem~6.9 in \cite{bsp}, the idea
  being as follows. First, the defining relations $abc=cba$ express the fact
  that the following operator must intertwine $u^{\otimes 3}$:
  $$T(e_i\otimes e_j\otimes e_k)=e_k\otimes e_j\otimes e_i$$

  The point is that $T$ comes from the following Brauer diagram:

  \setlength{\unitlength}{0.7pt}
  \begin{displaymath}
    p_3
    =\begin{matrix}
      \begin{picture}(40,60)
        \put(0,0){\line(1,3){5}}
        \put(5,15){\line(1,1){45}}
        \put(0,60){\line(1,-3){5}}
        \put(5,45){\line(1,-1){45}}
        \put(20,0){\line(3,4){22.5}}
        \put(20,60){\line(3,-4){22.5}}
      \end{picture}\end{matrix}
  \end{displaymath}
  
  A careful examination shows that this diagram ``generates'' all the diagrams
  having an even number of crossings, and this proves the result.
\end{proof}

In was pointed out in \cite{bsp} that $A_o^*(n)$ should appear as some kind of
``twist'' of $C(U_n)$, due to a certain common occurrence of the symmetrized
Rayleigh law, in the asymptotic representation theory of these algebras.  In
this paper we will present several results in this sense. These results will be
all based on the following fact.

\begin{theorem}
  If $u$, $v$ are the fundamental corepresentations of $A_o^*(n)$, $C(U_n)$ then
  $$\Hom(u^{\otimes k}, u^{\otimes l})=\Hom(v_k, v_l)$$
  for any $k$, $l$, where $v_k=v\otimes\bar v\otimes v\otimes\ldots$ ($k$ terms).
\end{theorem}

\begin{proof}
  If $\alpha$, $\beta$ are tensor products between $v$, $\bar{v}$, of length
  $K$, $L$, and we denote by $P(\alpha, \beta)\subset P(K, L)$ the set of
  pairings such that each string joins a $v$ to a $\bar{v}$, then:
  $$\Hom(\alpha, \beta)=\Span(T_p~|~ p\in P(\alpha, \beta))$$

  This is indeed a well-known result, see e.g. Theorem~9.1 of \cite{bco} for a
  recent proof of the version with non-crossing pairings. Now in the particular
  case $\alpha=v_k$, $\beta=v_l$, we get:
  $$\Hom(v_k, v_l)=\Span(T_p~|~ p\in P(v_k, v_l))$$

  On the other hand, Proposition~2.2 shows that we have $P(v_k, v_l)=E(k,
  l)$. Together with the middle assertion in Theorem~2.4, this gives the result.
\end{proof}

\section{The projective version}

The projective version of a unitary Hopf algebra $(A, u)$ is the subalgebra
$PA\subset A$ generated by the entries of $u\otimes\bar{u}$, with
$u\otimes\bar{u}$ as fundamental corepresentation.

In the Proposition below are some basic examples, with (1) justifying the
terminology. Recall that $A = C^*(\Gamma)$, with $\Gamma$ finitely generated, is
a Woronowicz algebra with fundamental corepresentation $u = \diag(g_1, \ldots,
g_n)$ where the $g_i$ are generators of $\Gamma$. On the other hand $A = C(G)$,
with $G \subset U_n$ compact subgroup, is a Woronowicz algebra with fundamental
corepresentation $u$ given by the embedding into $U_n$. We refer finally to
\cite{wa2} for the definition of the quantum automorphism group
$A_{aut}(M_n(\C))$.

\begin{proposition}
  The projective version is as follows:
  \begin{enumerate}
  \item For $G\subset U_n$ we have $PC(G)=C(PG)$.
  \item For $\Gamma=\langle g_i\rangle $ we have $PC^*(\Gamma)=C^*(\Lambda)$,
    with $\Lambda=\langle g_{i}g_j^{-1}\rangle \subset \Gamma$.
  \item For $A=A_o(n)$ and $A_u(n)$ we have $PA=A_{aut}(M_n(\C))$.
  \end{enumerate}
\end{proposition}

\begin{proof}
  The first two assertions are well-known, and follow from definitions. The
  third assertion is known as well: for $A_o(n)$, see Corollary~4.1 in
  \cite{ba3} and for $A_u(n)$, use Th\'eor\`eme~1~(iv) in \cite{ba1}.
\end{proof}

In this section we compute the projective version of $A_o^*(n)$. Our starting
point is the following simple observation, coming from definitions.
 
\begin{proposition}
  $PA_o^*(n)$ is commutative.
\end{proposition}

\begin{proof}
  This follows indeed from the relations $abc=cba$, because $PA_o^*(n)$ is
  generated by the elements of type $ab$, and we have $aba'b'=ab'a'b=a'b'ab$.
\end{proof}

\begin{theorem}
  We have $PA_o^*(n)=C(PU_n)$.
\end{theorem}

\begin{proof}
  We will prove that the algebras $A_o^*(n)$ and $C(U_n)$ have the same
  projective version. In order to ``compare'' these two algebras, we use Wang's
  universal algebra $A_u(n)$, having both of them as quotients. Consider indeed
  the following diagram:
  $$\xymatrix{A_u(n) \ar[d] \ar[r] & A^*_o(n) \ar[d] \\
    C(U_n) \ar[r] & C(O_n)} $$

  We fix $k$, $l\geq 0$ and we consider the words $\alpha =
  (u\otimes\bar{u})^{\otimes k}$ and $\beta = (u\otimes\bar{u})^{\otimes
    l}$. According to the above results, the spaces $\Hom(\alpha, \beta)$ of our
  four algebras appear as span of the operators $T_p$, with $p$ belonging to the
  following four sets of diagrams:
  $$\begin{matrix}
    N(2k, 2l)&\subset&E(2k, 2l)\\ \cap&&\cap\\
    E(2k, 2l)&\subset&P(2k, 2l)
  \end{matrix}$$

  Summarizing, we have computed the relevant diagrams for the projective
  versions of our four algebras. So, let us look now at these projective
  versions:
  $$\xymatrix{PA_u(n) \ar[d] \ar[r] & PA^*_o(n) \ar[d] \\
    C(PU_n) \ar[r] & C(PO_n)} $$

  We can see that the relationship between $PA_o^*(n)$ and $C(PU_n)$ is as
  follows:
  \begin{enumerate}
  \item These two algebras appear as quotients of a same algebra.
  \item The relevant diagrams for these two algebras are the same.
  \end{enumerate}

  It is a well-known application of Woronowicz's results in \cite{wo3} that
  these two conditions ensure the fact that our algebras are isomorphic and we
  are done.
\end{proof}

\section{A uniqueness result}

In this section we find an abstract characterization of the algebra $A_o^*(n)$:
it is in some sense the unique intermediate quotient $A_o(n) \to A \to
C(O_n)$. This will justify the terminology ``half-liberated'' that we use in
this paper.

In what follows, we use Woronowicz's tensor category formalism in
\cite{wo3}. That is, we call ``tensor category'' a tensor $C^*$-category with
duals whose monoid of objects is $(\N, +)$, embedded into the tensor
$C^*$-category of Hilbert spaces.  The following result follows from
Woronowicz's Tannakian duality in \cite{wo3}:

\begin{theorem}
  The intermediate Hopf algebras $A_o(n)\to A\to C(O_n)$ are in one-to-one
  correspondence with the tensor categories $C$ satisfying
  $$\Span(T_p~|~p\in N(k, l))\subset C(k, l)\subset \Span(T_p~|~p\in P(k, l))$$
  where $N$ denotes as usual the noncrossing pairings, and $P$, all the pairings.
\end{theorem}

We have then the following definition, adapted from \cite{bsp}.

\begin{definition}
  An intermediate Hopf algebra $A_o(n)\to A\to C(O_n)$ is called ``easy'' when
  its associated tensor category is of the form
  $$C(k, l)=\Span(T_p~|~ p\in D(k, l))$$
  for a certain collection of subsets $D(k, l)\subset P(k, l)$, with $k$, $l\in\N$.
\end{definition}

In other words, we know from Theorem~4.1 that the $\Hom$-spaces $C(k, l)$
associated to $A$ consist of certain linear combinations of partitions. In the
case where for any $k$, $l$ we can exhibit a certain set of partitions $D(k,
l)\subset P(k, l)$ such that $C(k, l)$ is spanned by the elements $T_p$ with
$p\in D(k, l)$, we call our Hopf algebra ``easy''.

Observe that the sets to be exhibited can be chosen to be:
$$D(k, l)=\{p\in P(k, l)~|~T_p\in C(k, l)\}$$

Thus, in concrete situations, the check of easiness is in fact
straightforward. We refer to \cite{bsp} for full details regarding this notion,
including examples and counterexamples. We use as well the following technical
definition.

\begin{definition}
  Let $p\in P(k, l)$ be a partition, with its points counted as usual
  counterclockwise starting from bottom left. For $i=1, 2, \ldots, k+l$ we denote
  by $p^i$ the partition obtained by connecting with a semicircle the $i$-th and
  $(i+1)$-th points.
\end{definition}

In this definition we agree of course that the points are counted modulo
$k+l$. The partitions $p^i$ will be called ``cappings'' of $p$, and will be
generically denoted $p'$. 

We denote by $P$, $E$, $N$ the collection of sets in Definition~2.1, endowed
with the operations of horizontal and vertical concatenation, and upside-down
turning. These operations are well-known to correspond via $p\to T_p$ to the
tensor product, composition and involution operations in the corresponding
$\Hom$-spaces: see e.g. Prop.~1.9 in \cite{bsp}.

\begin{lemma}
  Consider a partition $p\in P-N$, having $s\geq 4$ strings.
  \begin{enumerate}
  \item If $p\in P-E$, there exists a capping $p'\in P-E$.
  \item If $p\in E-N$, there exists a capping $p'\in E-N$.
  \end{enumerate}
\end{lemma}

\begin{proof}
  First, we can use a rotation --- see the proof of Lemma~2.7 in \cite{bsp} ---
  in order to assume that $p$ has no upper points. In other words, our data is a
  partition $p\in P(0, 2s)-N(0, 2s)$, with $s\geq 4$.

  (1) The assumption $p\notin E$ means that $p$ has certain strings having an
  odd number of crossings. We fix such an ``odd'' string, and we try to cap $p$,
  as for this string to remain odd in the resulting partition $p'$. An
  examination of all possible pictures shows that this is possible, provided
  that our partition has $s\geq 3$ strings, and we are done.

  (2) The assumption $p\notin N$ means that $p$ has certain crossing strings. We
  fix such a pair of crossing strings, and we try to cap $p$, as for these
  strings to remain crossing in $p'$. Once again, an examination of all possible
  pictures shows that this is possible, provided that our partition has $s\geq
  4$ strings, and we are done.
\end{proof}

For $p\in P$ we denote by $\langle p\rangle \subset P$ the collection of
partitions generated by $p$ and by $N$, via the above operations of
concatenation and upside-down turning. In particular if $q \in \langle p\rangle$
we also have $q' \in \langle p \rangle$ for any capping $q'$ of $q$. Observe
that we have:
$$\langle T_p\rangle =\Span(T_q~|~q\in \langle p\rangle )$$

Here the left term is by definition the tensor category generated by $T_p$. It
contains at least all morphisms $T_q$ for $q \in N$ since $T_|$ is the identity
morphism and $T_\sqcap$, $T_\sqcup$ are the morphisms describing the duality in
the category.

Let us quote two examples used for the proof of the next Lemma: for $p =
\big\slash\hskip-1.35ex\big\backslash$ we clearly have $\langle p \rangle = P$,
and if $p = p_3$ is the diagram pictured in the proof of Theorem~2.4, we have
$\langle p \rangle = E$ as already stated there.

\begin{lemma}
  Consider a partition $p\in P(k, l)-N(k, l)$.
  \begin{enumerate}
  \item If $p\in P(k, l)-E(k, l)$ then $\langle p\rangle =P$.
  \item If $p\in E(k, l)-N(k, l)$ then $\langle p\rangle =E$.
  \end{enumerate}
\end{lemma}

\begin{proof}
  This can be proved by recurrence on the number of strings,
  $s=(k+l)/2$. Indeed, by using Lemma~4.4, for $s\geq 4$ we have a descent
  procedure $s\to s-1$, and this leads to the situation $s\in\{1, 2, 3\}$, where
  the statement is clear from the examples above.
\end{proof}

\begin{theorem}
  $A_o^*(n)$ is the unique easy Hopf algebra between $A_o(n)$ and $C(O_n)$.
\end{theorem}

\begin{proof}
  Let $A$ be such an easy Hopf algebra, and consider the sets $D(k, l)\subset
  P(k, l)$, as in Definition~4.2. We have three cases:

  (1) Assume first that we have $D(k, l) = N(k, l)$, for any $k$, $l$. We can
  apply Theorem~4.1 and we get $A=A_o(n)$.

  (2) Now, assume that we have $D(k, l)\subset E(k, l)$ for any $k$, $l$, and
  that there exist $k'$, $l'$ and $p\in D(k', l')-N(k', l')$. From Lemma~4.5 (2)
  we get $\langle p\rangle =E$, and by applying Theorem~4.1 we get $A=A_o^*(n)$.

  (3) Assume finally that there exist $k$, $l$ and $p\in D(k, l)-E(k, l)$. From
  Lemma~4.5 (1) we get $\langle p\rangle =P$, and by applying Theorem~4.1 we get
  $A=C(O_n)$.
\end{proof}

\section{Diagonal groups}

In this section and in the next two ones we present a classification result for
the irreducible corepresentations of $A_o^*(n)$, which is reminiscent of the
classification by highest weights of the irreducible representations of compact
Lie groups.

We begin with some general considerations. In order to simplify the
presentation, all the Woronowicz algebras to be considered will be assumed to be
full.

\begin{theorem}
  Let $(A, u)$ be a Woronowicz algebra. Put
  $$A' = A \,/\, \left\langle u_{ij}=0, \forall i\neq j\right\rangle$$
  and denote by $g_i$ the image of $u_{ii}$ in $A'$.
  \begin{enumerate}
  \item $A'$ is a cocommutative Hopf algebra quotient and the unitaries $g_i$
    generate a group $L$ such that $A' \simeq C^*(L)$.
  \item If the elements $g_i \in L$ are pairwise distinct, then $A'$ is maximal
    as a cocommutative Hopf algebra quotient of $A$.
  \end{enumerate}
\end{theorem}

\begin{proof}
  (1) Denote by $J$ the closed, two-sided ideal generated by the $u_{ij}$ with
  $i\neq j$. Denote by $q : A \to A'$ the quotient map. We first have to prove
  that $(q\otimes q) \Delta (a) = 0$ for all $a\in J$, and it suffices to
  consider $a = u_{ij}$ with $i\neq j$. But then for any $k$ at least one of
  $q(u_{ik})$, $q(u_{kj})$ vanishes so that
  $$(q\otimes q)\Delta(u_{ij}) = \sum q(u_{ik})\otimes q(u_{kj}) = 0$$
  Hence $\Delta$ factors to a coproduct $\Delta' : A' \to A'\otimes
  A'$. Moreover the elements $g_i$ are group-like in $A'$:
  $$(q\otimes q)\Delta(u_{ii}) = \sum q(u_{ik})\otimes q(u_{ki}) = g_i\otimes g_i$$
  Since the $g_i$ generate $A'$, this shows that $A'$ is cocommutative.

  (2) Assume $q$ factors through another Hopf algebra quotient map $r : A
  \to A''$ with $A''$ cocommutative:
  \begin{displaymath}
    \xymatrix{A \ar[rd]_{q} \ar[r]^{r} & A'' \ar[d]^{s} \\ & A'}
  \end{displaymath}

  Denote by $u'$, $u''$ the images of $u$ in $M_n(A')$, $M_n(A'')$. We have by
  definition $u' = \diag(g_i)$, and since $A''$ is cocommutative $u''$ can
  also be decomposed into one-dimensional corepresentations: we write $u'' =
  P^{-1} \diag(h_i) P$ with $h_i \in A''$ and $P\in U_n$.

  By commutativity of the diagram above we have $(id\otimes s)(u'') = u'$ hence
  the elements $s(h_i)$ give the decomposition of $u'$ into irreducible
  subcorepresentations, so that we can find $\sigma\in S_n$ such that $s(h_i) =
  g_{\sigma(i)}$. Let $P_\sigma \in U_n$ denote the corresponding
  permutation matrix. 

  We have by construction
  \begin{eqnarray*}
    u' &=& (id\otimes s)(P^{-1} \diag(h_i) P) = P^{-1} \diag(s(h_i)) P \\
    &=&   P^{-1}P_\sigma^{-1} \diag(g_i) P_\sigma P =  P^{-1}P_\sigma^{-1} u' P_\sigma P
  \end{eqnarray*}
  Denoting $Q = P_\sigma P$ this yields $Q_{ij} g_j = g_i Q_{ij}$, and if the
  $g_i$ are pairwise distinct we obtain $Q_{ij} = 0$ for $i\neq j$. Hence we
  have $P = P_\sigma$ up to scalar factors, and the identity $u'' = P^{-1}
  \diag(h_i) P$ shows that $u''$ was already diagonal. As a result $\Ker r
  \supset J$, so we have in fact equality and $s$ is an isomorphism.
\end{proof}

\begin{definition}
  The discrete group $L$ given by
  $$C^*(L)=A\,/\,\langle u_{ij}=0,\,\forall i\neq j\rangle$$
  is called diagonal group of the Woronowicz algebra $A$.
\end{definition}

Let us now discuss the diagonal groups for standard examples.

\begin{proposition}
  We have the following results:
  \begin{enumerate}
  \item For $A=C^*(\Gamma)$ with $\Gamma$ finitely generated discrete group we
    have $L=\Gamma$.
  \item For $A=C(G)$ with $G\subset U_n$ compact we have $L=\widehat{T}$,
    where $T=G\cap\T^n$.
  \item For $A=A_u(n)$, $A_o(n)$, $A_s(n)$ we have $L=F_n$, $\Z_2^{*n}$, $\{1\}$
    respectively.
  \end{enumerate}
\end{proposition}

\begin{proof}
  This is clear from definitions.
\end{proof}

The interest in the diagonal group comes from Proposition~5.3 (2): for $C(U_n)$,
this group is nothing but the dual of the maximal torus of $U_n$. For a general
connected compact subgroup $G \subset U_n$, the diagonal group need not be a
dual maximal torus, e.g. $O_n \cap \T^n = \Z_2^n \subset U_n$ is maximal abelian
but not a torus. However all maximal torii are clearly duals of diagonal groups,
up to conjugation of the fundamental representation $u$ by a matrix $P \in U_n$,
and they are known to be maximal abelian. For $G = S_n \subset U_n$ we have $G
\cap \T^n = \{1\}$, which is not maximal abelian.

In the quantum case there is no clear notion of what a torus should be, however
there are cases where diagonal groups are clearly too small, e.g. for $A =
A_s(n)$. We will attack the issue by considering the potential applications of
diagonal groups to representation theory: we introduce below a map $\Phi$ which
should be injective for ``good'' diagonal groups.

We denote by $R^+(A)$ the set of equivalence classes of finite dimensional
smooth corepresentations of $A$, endowed with the operations of sum and tensor
product. We use the character map $\chi : R^+(A)\to A$, given by
$\chi(r)=Tr(r)$. 

\begin{definition}
  Associated to a Woronowicz algebra $(A, u)$ is the map
  $$\Phi:R^+(A)\to\N[L]$$
  given by $r\to\chi(r)'$, where $L$ is the diagonal group. 
\end{definition}

In this definition $x\to x'$ is the canonical map $A\to A'=C^*(L)$, constructed
in the previous section. Observe that an alternative definition for $\Phi$ could
be $\Phi(r)=\chi(r')$, where $r'\in M_n(A')$ is the corepresentation induced by
$r\in M_n(A)$.

The target of $\Phi$ is indeed $\N[L]$, because characters of corepresentations
of $C^*(L)$ are sums of elements of $L$. For the same reason the elements
$\Phi(r) \in \N[L]$ can also be considered as subsets with repetitions of $L$,
which we will denote by $\Sigma(r)$.

Observe finally that $\Phi$ is a morphism of semirings, due to the additivity
and multiplicativity properties of the character map $w\to\chi(w)$.

\begin{theorem}
  In the following situations, $\Phi$ is injective and $C^*(L)$ is a maximal
  cocommutative quotient:
  \begin{enumerate}
  \item For $A=C^*(\Gamma)$, with $\Gamma$ discrete group of finite type.
  \item For $A=C(G)$, with $G \subset U_n$ connected such that $G \cap
    \T$ is a maximal torus.
  \item For the free quantum algebras $A_o(n)$, $A_u(n)$.
  \end{enumerate}
\end{theorem}

\begin{proof}
  We use the explicit computations of $L$ given at Proposition~5.3.
  \begin{enumerate}
  \item Here $L=\Gamma$, and the quotient map as well as $\Phi$ are actually
    isomorphisms.
  \item Here $L$ is the weight lattice of $G$ with respect to the maximal torus
    $T = G \cap \T$, and $\Phi$ is the character map, which is known to classify
    representations of $G$. The maximality result holds because maximal torii
    are maximal abelian.
  \item Here $L=\Z_2^{*n}$, $F_n$, and the injectivity is easily proved using
    the fusion rules of $A_o(n)$, $A_u(n)$. The maximality results from
    Theorem~5.1~(2).
  \end{enumerate}
\end{proof}

\section{Representation theory}

We have seen in the previous section that, at least for certain Woronowicz
algebras and up to conjugation of the fundamental corepresentation, the diagonal
group is a reasonable candidate for a ``dual maximal torus''. More precisely, we
can say that we have a dual maximal torus when the quotient is maximal
cocommutative and $\Phi$ is injective.

In what follows we will prove that these requirements are fulfilled in the case
of the algebra $A_o^*(n)$. This result, besides of being of independent
theoretical interest, can be regarded as a concrete classification of the
corepresentations of $A_o^*(n)$, in terms of ``combinatorial data''.

We begin with a study of the diagonal group. The next Proposition shows in
particular that Theorem~5.1 (2) applies in the case of $A_o^*(n)$, hence the
diagonal quotient $C^*(L_n)$ of $A_o^*(n)$ is maximal cocommutative. Recall that
the diagonal group $L_n$ of $A_o^*(n)$ was already introduced at Definition~1.4.

\begin{proposition}
  Write $\Z_2=\{1, \tau\}$ and let $\tau$ act on $\Z^n$ by $\tau\cdot\lambda =
  -\lambda$. Consider the subsets $L_n^\circ = \{(\lambda_i)\cdot
  1~|~\sum\lambda_i=0\}$ and $L_n^\tau =
  \{(\lambda_i)\cdot\tau~|~\sum\lambda_i=1\}$ of $\Z^n\rtimes \Z_2$.
  \begin{enumerate}
  \item The group $L_n$ embeds into $\Z^n\rtimes\Z_2$, via
    $g_i=e_i\cdot\tau$.
  \item Its image is $L_n^\circ \cup L_n^\tau$.
  \item We have $L_n \simeq \Z^{n-1} \rtimes \Z$.
  \end{enumerate}
\end{proposition}

\begin{proof}
  It follows from the definition of the semidirect product that the elements
  $\gamma_i=e_i\cdot\tau$ multiply according to the following formula:
  $$\gamma_{i_1}\ldots\gamma_{i_k}=
  \left(e_{i_1}-e_{i_2}+\ldots+(-1)^{k+1}e_{i_k}\right)\cdot\tau^k$$ 
  In particular with $k=2$, $3$ we get: $\gamma_a\gamma_b = (e_a-e_b)\cdot 1$,
  $\gamma_a\gamma_b\gamma_c = (e_a-e_b+e_c)\cdot\tau$.  Thus we can define a
  morphism $\varphi:L_n\to \Z^n\rtimes\Z_2$ by
  $\varphi(g_i)=\gamma_i$. Moreover, the above formula shows that the image of
  $\varphi$ is the subgroup in the statement.

  If $w$ is a word on $g_1, \ldots, g_n$, we denote by $w^{\rm odd}$, $w^{\rm
    even}$ the subwords formed by letters at odd and even positions
  respectively, and by $w_i$ the number of occurrences of $g_i$ in $w$. Then the
  map $\varphi$ is given by the following formula, where $x$ is the unique
  element of $\Z_2$ making $\varphi(w)$ an element of $\varphi(L_n)$:
  $$\varphi(w)=(w^{\rm odd}_i-w^{\rm even}_i)_i\cdot x$$  
  Indeed, the above formula holds for $w=g_i$, and an easy computation shows
  that the expression on the right is multiplicative in $w$.

  Now we can prove that $\Ker(\varphi)$ is trivial. If a word $w$ lies in
  $\Ker(\varphi)$, the above formula shows that each $g_i$ appears an equal
  number of times at odd and even positions of $w$. But by definition of $L_n$
  the letters of $w^{\rm odd}$ can be permuted without changing the group element,
  hence we can bring pairs of $g_i$'s in $w^{\rm even}$ and $w^{\rm odd}$ side-by-side
  and simplify them according to the relation $g_i^2=1$, and we get $w=1$.

  Finally it is clear that $\varphi(L_n) = \{xy~|~x\in L_n^\circ,\, y\in \{1,
  e_1\cdot\tau\} \}$. Since $L_n^\circ \simeq \Z^{n-1}$ and $\{1, e_1\cdot\tau\}
  \simeq \Z_2$ we have $\varphi(L_n) \simeq \Z^{n-1}\rtimes \Z_2$ and an easy
  check shows that the action of $\Z_2$ on $\Z^{n-1}$ is indeed given by $\tau
  \cdot \lambda = -\lambda$.
\end{proof}

In this section and the next ones we will make frequent use of the following
map, which connects, in a sense to be precised, the corepresentation theory of
$A_o^*(n)$ to the representation theory of $U_n$:
$$\psi : L_n \to \Z^n,~ (\lambda_i)\cdot x \mapsto (\lambda_i)$$
Note that $\psi$ is injective, and that it is not a group morphism.

\begin{theorem}
  $\Phi$ is injective for the algebra $A_o^*(n)$.
\end{theorem}

\begin{proof}
  Let $v$ be the fundamental corepresentation of $C(U_n)$, and consider the
  $k$-fold tensor product $v_k=v\otimes\bar v\otimes v\otimes\ldots$ According
  to Theorem~2.5, we have:
  $$\End(u^{\otimes k})=\End(v_k)$$

  Now recall that the subcorepresentations of a corepresentation $w$ are of the
  form $(p\otimes 1)w$, with $p\in \End(w)$ projection. This shows that we have
  a one-to-one additive correspondence $J$ between subobjects of $u^{\otimes k}$
  and subobjects of $v_k$, by setting:
  $$J((p\otimes 1)u^{\otimes k})=(p\otimes 1)v_k$$

  We claim that when $k$ varies, these $J$ maps are compatible with each
  other. Indeed, let $p$, $q$ be projections yielding irreducible
  subrepresentations of $v_k$, $v_l$. The same diagrammatic identifications as
  before show that:
  $$q\Hom(u^{\otimes l}, u^{\otimes k})p=q\Hom(v_l, v_k)p$$

  This shows that $(p\otimes 1)u^{\otimes k}=(q\otimes 1)u^{\otimes l}$ is
  equivalent to $(p\otimes 1)v_k=(q\otimes 1)v_l$, so the $J$ maps are indeed
  compatible with each other.  Summarizing, we have constructed an embedding of
  additive monoids:
  $$J:R^+(A_o^*(n))\to R^+(C(U_n))$$

  Now consider the map $\psi : L_n\to\Z^n$ introduced above, and extend it by
  linearity to $\N[L_n]$ and $\C[L_n]$. We claim that the following diagram is
  commutative, so that the injectivity of $\Phi$ for $A_o^*(n)$ follows from the
  one for $C(U_n)$, which is known from the classical theory:
  $$\xymatrix{  R^+(A_o^*(n)) \ar[d]^\Phi \ar[r]^J & R^+(C(U_n)) \ar[d]^\Phi \\
    \N[L_n] \ar[r]^\psi & \N[\Z^n]}$$

  Indeed, let us consider the linear extension $\psi : \C[L_n]\to
  \C[\Z^n]$. Since the quotient map $C(U_n) \to C^*(\Z^n)$ maps $v_{ij}$ to
  $\delta_{ij} e_i$, we have the following computation, with $\tilde{v}=v$ or
  $\bar{v}$ depending on the parity of $k$:
  \begin{eqnarray*}
    \psi((u_{i_1j_1}\ldots u_{i_kj_k})')
    &=&\delta_{i_1j_1}\ldots\delta_{i_kj_k}\psi(g_{i_1}\ldots g_{i_k})\\
    &=&\delta_{i_1j_1}\ldots\delta_{i_kj_k}(e_{i_1}-e_{i_2}+\ldots+(-1)^{k+1}e_{i_k})\\
    &=&(v_{i_1j_1}\bar{v}_{i_2j_2}\ldots\tilde{v}_{i_kj_k})'
  \end{eqnarray*}
  This shows that for any rank one projection $p$ we have:
  $$\psi(((p\otimes 1)u^{\otimes k})')=((p\otimes 1)v_k)'$$
  By linearity this formula must hold for any $p\in \End(u^{\otimes k})$, so we
  get $\psi(\Phi(r))=\Phi(J(r))$ for any corepresentation $r\subset u^{\otimes
    k}$, and the diagram commutes.
\end{proof}

\section{Highest weights}

We know from the previous section that the corepresentations of $A_o^*(n)$ can
be indexed by certain elements of $\N[L_n]$. In this section we furter
refine this result, by indexing the irreducible corepresentations of $A_o^*(n)$
by their ``highest weights''.

We first recall the general theory for $U_n$. With the choice of the basis
$(e_i-e_{i+1})_i$ for the root system associated to $T=\T^n$, the objects
of interest are as follows.

\begin{definition}
  Associated to $U_n$ are the following objects.
  \begin{enumerate}
  \item Dual maximal torus: $X=\Z^n$.
  \item Root system: $X_*=\{e_i-e_j ~|~ i\neq j\}$.
  \item Root lattice: $X^\circ=\{(\lambda_i)\in X ~|~ \sum\lambda_i=0\}$.
  \item Positive weights: $X_+=\{(\lambda_i)\in X^\circ~|~\lambda_1\geq
    0,\lambda_1+\lambda_2\geq 0,\ldots,\sum_1^{n-1}\lambda_i\geq
    0\}$.
  \item Dominant weights: $X_{++}=\{(\lambda_i)_i\in
    X~|~\lambda_1\geq\lambda_2\geq\ldots\geq\lambda_n\}$.
  \end{enumerate}
\end{definition}

Recall that $X_* \cup \{0\}$ is the set of weights $\Sigma(v\otimes \bar v)$ of
the adjoint representation of $U_n$, with the notation $\Sigma(\cdot)$
introduced after Definition~5.4. The set $X_*$ generates the root lattice
$X^\circ = \{(\lambda_i) ~|~ \sum\lambda_i=0\}\subset X$ and $X = \bigsqcup_{k\in\Z}
(X^\circ + k e_1)$.

Observe that $X_+$ is the subset of $X^\circ$ consisting of elements which have
positive coefficients with respect to $(e_i - e_{i+1})_i$. We endow $X$ with the
partial order $x\geq y$ if $x-y\in X_+$.

\begin{theorem}
  For any irreducible representation $w$ of $U_n$, the set with repetitions
  $\Sigma(w)$ has a greatest element $\lambda_w\in X$, called highest weight of
  $w$. For $w$, $w'$ irreducible we have $w\simeq w'$ iff
  $\lambda_w=\lambda_{w'}$. Finally, the set of highest weights is $X_{++}$.
\end{theorem}

\begin{proof}
  This is a well-known result concerning the classification of irreducible
  representations of $U_n$.
\end{proof}

We make the following definition.

\begin{definition}
  Associated to $A_o^*(n)$ are the following objects.
  \begin{enumerate}
  \item Dual maximal torus: $L=L_n\subset\Z^n\rtimes\Z_2$.
  \item Root system: $L_*=\{(e_i-e_j)\cdot 1 ~|~ i\neq j\}$.
  \item Root lattice: $L^\circ=\{(\lambda_i)_i\cdot 1\in L ~|~ \sum\lambda_i=0\}$.
  \item Positive weights: $L_+=\{(\lambda_i)_i\cdot 1\in L^\circ ~|~\lambda_1\geq
    0,\lambda_1+\lambda_2\geq 0,\ldots,\sum_1^{n-1}\lambda_i\geq 0\}$.
  \item Dominant weights: $L_{++}=\{(\lambda_i)_i\cdot x\in
    L~|~\lambda_1\geq\lambda_2\geq\ldots\geq\lambda_n\}$.
  \end{enumerate}
\end{definition}

Observe that we have as in the classical case $L_* \cup \{0\} = \Sigma(u\otimes
\bar u)$, where $u = \bar u$ is the fundamental corepresentation of $A_o^*(n)$.
The set $L_*$ generates the subgroup $L^\circ = L_n^\circ$ and we have $L =
L^\circ \sqcup L^\tau$ with $L^\tau = L_n^\tau = \{(\lambda_i)\cdot\tau ~|~ \sum
\lambda_i = 1\}$.

Note also that $L_+$ in contained in $L^\circ$ and consists of elements having
positive coefficients with respect to the basis $(e_i - e_{i+1})_i \subset L_*$
of $L^\circ$. We endow $L$ with the partial order $x\geq y$ if $xy^{-1}\in L_+$.

\begin{theorem}
  For any irreducible corepresentation $w$ of $A_o^*(n)$, the set with
  repetitions $\Sigma(w)$ has a greatest element $\lambda_w\in L$, called
  highest weight of $w$. For $w$, $w'$ irreducible we have $w\simeq w'$ iff
  $\lambda_w=\lambda_{w'}$. Finally, the set of highest weights is $L_{++}$.
\end{theorem}

\begin{proof}
  Recall that we use the map $\psi : L \to X$, $(\lambda_i)\cdot x \mapsto
  (\lambda_i)$ and observe that we have $\psi(L_+)=X_+$ and
  $$\psi((\lambda\cdot 1)(\mu\cdot x))=\psi(\lambda\cdot 1)+\psi(\mu\cdot x)$$
  Thus $\psi$ respects the orders. Now since $\Sigma(w)$ has a greatest element
  for any irreducible representation $w$ of $U_n$, and this highest weight
  characterizes $w$ up to isomorphism, the first two assertions follow as in
  Theorem~6.2.

  Moreover, we see that the images under $\psi$ of the highest weights of
  $A_o^*(n)$ are precisely the highest weights of the irreducible subobjects of
  the representations $v_k$ of $U_n$.

  Restricting representations of $U_n$ to $Z(U_n)=\T$ corresponds to projecting
  weights in the quotient $\Z=\Z^n/\langle e_i=e_{i+1}\rangle $, and hence we
  see that the weights $(\lambda_i)$ of $v^{\otimes k}\otimes\bar{v}^{\otimes
    l}$ satisfy $\sum\lambda_i=k-l$. Since all the irreducible representations
  appear as subobjects of some $v^{\otimes k}\otimes \bar{v}^{\otimes l}$, we
  conclude that the highest weights of the subobjects of $v_{2k}$
  (resp. $v_{2k+1}$) are exactly the dominant weights $(\lambda_i)\in X_{++}$
  such that $\sum\lambda_i=0$ (resp. $1$).

  Thus the highest weights we are looking for are the dominant weights of $U_n$
  which lie in the image of $\psi$, and the result is proved.
\end{proof}

As a first application, observe that the highest weights in $L_{++}^\circ$
correspond to the irreducible subobjects of even tensor powers $u^{\otimes 2k}$,
and they are mapped by $\psi$ to the highest weights of subobjects of tensor
powers of the adjoint representation $v\otimes\bar v$ of $U_n$. Hence we recover
the identification $PA_o^*(n)=C(PU_n)$ from Theorem~3.3, at the level of
irreducible corepresentations.

\section{Fusion rules}

In this section we explain how the fusion rules of $A_o^*(n)$ are related to
those of $U_n$, in terms of highest weights. From this we will be able to draw
the Cayley graph associated with $R^+(A_o^*(n))$. This is the point where the
group structure of the ``weight lattice'' $L$, and in particular its
non-commutativity, come into play.

According to Theorem~7.4, decomposing into irreducibles and taking highest
weights yields an additive bijection from $R^+(A_o^*(n))$ to the free abelian
semigroup $\N[L_{++}]$ generated by $L_{++}$. We endow $\N[L_{++}]$ with the
associative, biadditive product $\otimes$ given by the tensor product in
$R^+(A_o^*(n))$, and with the additive involution $\lambda \mapsto \bar\lambda$
given by the conjugation in $R^+(A_o^*(n))$. We proceed similarly for $U_n$ in
$\N[X_{++}]$.

\begin{theorem}
  Denote by $\psi : \N[L_{++}] \to \N[X_{++}]$ the natural additive injective
  map given by $\psi((\lambda_i)_i \cdot x)= (\lambda_i)_i$ on $L_{++}$. For
  $\nu\in X$, $\lambda \in L$ put $\nu^\lambda = \bar\nu$ if $\lambda \in
  L^\tau$ and $\nu^\lambda = \nu$ if $\lambda \in L^\circ$. Then we have for all
  $\lambda$, $\mu \in L_{++}$ the following equalities in $\N[X_{++}]$
  \begin{eqnarray*}
    \psi(\lambda\otimes\mu) &=& \psi(\lambda)\otimes \psi(\mu)^\lambda \\
    \psi(\bar\lambda) &=& \overline{\psi(\lambda)}^\lambda
  \end{eqnarray*}
\end{theorem}

\begin{proof}
  Denote by $v_\lambda$ (resp. $w_{\psi(\lambda)}$) an irreducible corepresentation of
  $A_o^*(n)$ (resp. representation of $U_n$) with highest weight $\lambda \in
  L_{++}$ (resp. $\psi(\lambda) \in X_{++})$. 

  It is clear from the definition of the maps $\Phi$ that we have
  \begin{eqnarray*}
    \Phi(v_\lambda\otimes v_\mu) &=& \Phi(v_\lambda)\Phi(v_\mu) \in \N[L] \\
    \Phi(w_{\psi(\lambda)}\otimes w_{\psi(\mu)}) &=& 
    \Phi(w_{\psi(\lambda)})\Phi(w_{\psi(\mu)}) \in \N[X] \\
    \Phi(w_{\psi(\mu)}^\lambda) &=& \Phi(w_{\psi(\mu)})^\lambda
  \end{eqnarray*}
  where $w^\lambda = \bar w$ or $w$ and $\nu^\lambda = \pm \nu \in X$ are the
  natural actions of $L$ on $R^+(U_n)$ and $X$ via $\Z_2$. Moreover we
  know from the proof of Theorem~6.2 that $\psi(\Phi(v_\lambda)) =
  \Phi(w_{\psi(\lambda)})$ for all $\lambda \in X_{++}$. 

  Hence by injectivity of $\Phi$ it suffices to show that
  \begin{displaymath}
    \psi(\Phi(v_\lambda)\Phi(v_\mu)) = \psi(\Phi(v_\lambda)) \psi(\Phi(v_\mu))^\lambda
  \end{displaymath}
  
  We observe that $\Phi(v_\lambda) \subset \N [L^x]$ if $\lambda \in
  L_{++}^x$, for $x = 1$, $\tau$: this is e.g. a consequence of the fact that
  the elements of $\Sigma(v_\lambda)$ are smaller than $\lambda$. Now the result
  follows from the identity $\psi(\nu \nu') = \psi(\nu)\psi(\nu')^\nu$ which is
  clear from the definition of the product in $L \subset \Z^n \rtimes
  \Z_2$.

  Denote by $(x\to x^\sharp)$ the linear extension of $(\lambda \mapsto
  \lambda^{-1})$ to $\N[L]$, and of $(\lambda\to -\lambda)$ to $\N[X]$. We have
  then
  \begin{eqnarray*}
    \Phi(\bar v_\lambda) &=& \Phi(v_\lambda)^\sharp \\
    \Phi(\bar w_{\psi(\lambda)}) &=& \Phi(w_{\psi(\lambda)})^\sharp 
  \end{eqnarray*}
  Hence the second identity of the statement follows as the first one from the
  identity $\psi(\nu^\sharp) = (\psi(\nu)^\sharp)^\nu$.
\end{proof}

We see in particular that all corepresentations of $A_o^*(n)$ with weights in
$L^\tau$ are selfadjoint. As an other application let us prove that the fusion
rules of $A_o^*(n)$ are not commutative. This might seem quite surprising,
because $A_o^*(n)$ appears as intermediate subalgebra between $A_o(n)$ and
$C(O_n)$, both having commutative fusion rules.

\begin{proposition}
  For $n\geq 3$ the fusion rules of $A_o^*(n)$ are not commutative. In the case
  $n=3$, we have $u\otimes w \not\simeq w\otimes u$, where $u$ is the
  fundamental corepresentation, and $w$ is the irreducible subobject of
  $u^{\otimes 4}$ with highest weight $(1,1,-2)\cdot 1$.
\end{proposition}

\begin{proof}
  First recall that
  \begin{eqnarray*}
    \psi(L^\circ_{++}) &=&
    \{(\lambda_i)_i \in X_{++} ~|~ \textstyle\sum \lambda_i = 0\} \\
    \psi(L^\tau_{++}) &=&
    \{(\lambda_i)_i \in X_{++} ~|~ \textstyle\sum \lambda_i = 1\}
  \end{eqnarray*}
  Thus by Theorem~8.1 it is enough to produce irreducible representations $v$,
  $w$ of $U_n$ with highest weights in $\psi(L^\tau_{++})$, $\psi(L^\circ_{++})$
  respectively such that $v\otimes w \not\simeq v\otimes\bar w$. 

  In fact for any irreducible representations $v$, $w$, $w'$ of $U_n$ we have
  $v\otimes w \simeq v\otimes w'$ iff $w \simeq w'$. Indeed if $\lambda$, $\mu$
  are the highest weights of $v$, $w$, the (non-irreducible) representation
  $v\otimes w$ admits $\lambda+\mu$ as highest weight.

  In particular if $u$ is the fundamental corepresentation of $A_o^*(n)$, and
  $w$ is an irreducible corepresentation of $C(PU_n) \simeq PA_o^*(n) \subset
  A_o^*(n)$ such that $\bar w \not\simeq w$, we have $u\otimes w \not\simeq
  w\otimes u$ in $R^+(A_o^*(n))$. Such representations $w$ of $PU_n$ always
  exist when $n\geq 3$. 

  In the case $n=3$, the tensor square $u^{\otimes 2}$ is the sum of the trivial
  representation and a selfadjoint representation, but the power $u^{\otimes 4}$
  has nonselfadjoint irreducible subrepresentations. Recall indeed that
  $\lambda_{\bar w} = -w_0(\lambda_w)$, where $w_0(\lambda_1,\ldots,\lambda_n) =
  (\lambda_n,\ldots,\lambda_1)$.  Hence $w_1$, $w_2$ with highest weights
  $(1,1,-2)$, $(2,-1,-1)$ are subrepresentations of $u^{\otimes 4}$ which are
  nonselfadjoint.
\end{proof}

\bigskip

Now let us turn to Cayley graphs. Recall the following definition:
\begin{definition}
  Let $(A, u)$ be a Woronowicz algebra and fix a self-adjoint corepresentation
  $u_1$ of $A$ not containing the trivial corepresentation. The Cayley graph
  associated to this data is defined as follows:
  \begin{enumerate}
  \item the vertices are classes of irreducible corepresentations,
  \item if $w\subset v\otimes u_1$, we draw $\dim \Hom(w, v\otimes u_1)$ edges
    from $v$ to $w$.
  \end{enumerate}
\end{definition}

In the case of $C(U_n)$ we take $u_1 = v \oplus \bar v$, where $v$ is the
fundamental representation of $U_n$ on $\C^n$. In the case of $C(PU_n)$
we take for $u_1$ the unique nontrivial irreducible subrepresentation of
$v\otimes \bar v$. In the case of $A_o^*(n)$ we take $u_1 = u = \bar u$, where
$u$ is the fundamental corepresentation. Identifying irreducible
corepresentations with their highest weights we can draw the corresponding
Cayley graphs in $L_{++}$, $X_{++}$. We moreover identify $L_{++}$ with a subset
of $X_{++}$ via the map $\psi$ as previously.

We refer to the pictures of Section~9 for an illustration of the Propositions
below in the case $n=3$.

\begin{proposition}
  The Cayley graph of $A_o^*(n)$ corresponds to the full subgraph of the Cayley
  graph of $C(U_n)$ whose vertices are elements $(\lambda_i)_i$ lying in
  $\psi(L_{++})$, i.e. such that $\sum \lambda_i = 0$ or $1$. In particular
  there is an edge of the Cayley graph of $A_o^*(n)$ between two such elements
  $\lambda$, $\mu$ of $X_{++}$ {\em iff} $\lambda - \mu = \pm e_i$ for some $i$.
\end{proposition}

\begin{proof}
  We already know that the vertices of the Cayley graph of $A_o^*(n)$ correspond
  to the image of $\psi : L_{++} \to X_{++}$, i.e. to the elements
  $(\lambda_i)_i \in X_{++}$ such that $\sum \lambda_i \in\{0,1\}$.

  Let $(\lambda, \mu)$ be an edge between two such vertices in the Cayley graph
  of $C(U_n)$. Let $w_\lambda$, $w_\mu$ be irreducible representations of $U_n$
  with highest weights $\lambda$, $\mu$. We have $w_\mu \subset w_\lambda
  \otimes v$ or $w_\mu \subset w_\lambda \otimes \bar v$, hence $\mu \in
  \Sigma(w_\lambda) + \{e_i\}$ or $\mu \in \Sigma(w_\lambda) + \{-e_i\}$
  respectively, which yields $\sum \mu_i = \sum \lambda_i \pm 1$ because $\sum
  \nu_i$ is constant for $\nu\in\Sigma(w_\lambda)$.

  So we have $w_\mu \subset w_\lambda\otimes v$ if $\lambda \in \psi(L_{++}^\circ)$
  and $w_\mu \subset w_\lambda\otimes \bar v$ otherwise. According to
  Theorem~8.1 this reads $w'_\mu \subset w'_\lambda \otimes u$ in both cases,
  where $w'_\lambda$, $w'_\mu$ are irreducible corepresentations of $A_o^*(n)$
  with highest weights $\lambda$, $\mu$. Hence $(\lambda, \mu)$ is also an edge
  in the Cayley graph of $A_o^*(n)$.

  The last assertion holds for $A_o^*(n)$ because it holds for $U_n$. As a
  matter of fact if $w$ is a representation of $U_n$ with highest weight
  $\lambda \in X_{++}$, then the highest weights of irreducible subobjects of
  $w\otimes (v\oplus \bar v)$ are the elements of $\{\lambda \pm e_i\} \cap
  X_{++}$.
\end{proof}

\begin{proposition}
  Consider in the Cayley graph of $A_o^*(n)$ the vertices with highest weight
  $\lambda \in L_{++}^\circ$, and the paths of length $2$ between such
  vertices. Remove one loop at each vertex. Then the graph obtained coincides
  with the Cayley graph of $C(PU_n)$.
\end{proposition}

\begin{proof}
  This results clearly from the identification between $PA_o^*(n)$ and
  $PU_n$. Note that paths of length $2$ from $w$ to $w'$ correspond to
  inclusions $w' \subset w\otimes v\otimes\bar v$, and we have in $PU_n$ the
  decomposition $w\otimes v\otimes\bar v = w \oplus (w\otimes u_1)$, thus we
  obtain one more loop at each vertex with paths of length $2$ than in the
  Cayley graph of $C(PU_n)$ generated by $u_1$.
\end{proof}

\section{Polynomial growth}

Recall the following notion of growth introduced in \cite{ve2} and
\cite{bv1}. We fix a Woronowicz algebra $(A, u)$ and a self-adjoint
corepresentation $u_1$ not containing the trivial corepresentation. For any
irreducible corepresentation $w$ of $(A, u)$ we denote by $l(w)$ the length of
$w$, which is the distance from $1$ to $w$ in the Cayley graph associated with
$A$ and $u_1$. In other words we have
$$l(w) = \min \{k\in\N ~|~ w \subset u_1^{\otimes k}\}$$

Then the volumes of balls in the discrete quantum group associated with $(A, u)$
are defined as follows:
$$b_k = \sum_{l(w)\leq k} \dim(w)^2$$

\begin{definition}
  We say that $(A, u)$ has polynomial growth if the sequence $(b_k)$ has
  polynomial growth. If there exist constants $d$, $\alpha$, $\beta > 0$ such
  that $\alpha k^d \leq b_k \leq \beta k^d$ for all $k$, we say that $(A, u)$
  has polynomial growth with exponent $d$. These notions are independent of
  $u_1$.
\end{definition}

\bigskip

To prove the next Theorem we will work with representations and highest weights
of $SU_n$. Observe that restricting representations of $U_n$ to $SU_n$ amounts
at the level of weights to quotienting $X = \Z^n$ by the line generated by $e_1
+ \cdots + e_n$. We denote by $\bar X$ the quotient lattice and use it as the
lattice of weights of $SU_n$.

Since all irreducible representations of $SU_n$ are restrictions of irreducible
representations of $U_n$, the set of dominant weights $\bar X_{++} \subset \bar
X$ is the image of $X_{++}$ in $\bar X$. Moreover the quotient map is injective
on the subset of dominant weights $(\lambda_i)_i$ such that $0\leq \sum
\lambda_i \leq n-1$, and in particular it is injective on the image of $\psi :
L_{++} \to X_{++}$. Notice that the embedding
$$L_{++}^\circ \subset L_{++} \longrightarrow X_{++} \longrightarrow \bar X_{++}$$
is induced by the canonical embedding $R^+(C(PU_n)) \to R^+(C(SU_n))$ modulo the
identification $PA_o^*(n) \simeq C(PU_n)$.

Moreover the irreducible corepresentations of $C(SU_n)$, $C(U_n)$, $C(PU_n)$ and
$A_o^*(n)$ that correspond to each other in this picture have the same
dimensions: this is clear for the classical groups since restricting or
factoring a representation does not change its dimension, and for $A_o^*(n)$
this results from the proof of Theorem~6.2, where the correspondence between two
irreducible corepresentations $w$, $w'$ of $C(U_n)$, $A_o^*(n)$ is defined by the
common subspace $p(\C^{\otimes nk})$ on which they act.

\begin{theorem}
  $A_o^*(n)$ has exponential growth with exponent $d = n^2-1$.
\end{theorem}

\begin{proof}
  We proceed by comparison with $SU_n$, whose growth exponent $d = n^2-1$ is
  known by Theorem~2.1 in \cite{bv1}.
  
  Denote by $B_k^{PU}$, $B_k^{SU}$, $B_k^A$ the balls of radius $k$ in the
  Cayley graphs of $C(PU_n)$, $C(SU_n)$, $A_o^*(n)$, and by $b_k^{PU}$,
  $b_k^{SU}$, $b_k^A$ the growth sequences as above. We have $B_k^{PU} \subset
  B_{2k}^A \subset B_{2k}^{SU}$ by Propositions~8.4 and~8.5, hence $b_k^{PU}
  \leq b_{2k}^A \leq b_{2k}^{SU}$. 

  Now it remains to control $(b^{SU}_k)$ by $(b^{PU}_k)$, and this is a problem
  in the classical representation theory of $SU_n$. We present here an ad-hoc
  argument.

  Let first $w$ be an irreducible representation of $PU_n$, hence also of
  $SU_n$, such that $w \in B_k^{SU}$. Since the fundamental representation $v$
  of $SU_n$ generates its category of representations, we can find a constant
  $a$ depending only on $n$ and $l \leq ak$ such that $w\subset v^{\otimes
    l}$. Since $w$ is a representation of $PU_n$ we must have $l = pn$ for some
  $p\in\N$, and we notice that
  $$v^{\otimes n} = v^{\otimes n-1}\otimes v \subset 
  v^{\otimes n-1}\otimes \bar v^{\otimes n-1} = (v\otimes\bar v)^{\otimes n-1}$$
  and hence $w \subset (v\otimes\bar v)^{\otimes p(n-1)} \subset (v\otimes\bar
  v)^{\otimes l}$, so $w\in B_{ak}^{PU}$.

  Now take an irreducible representation $w\in B_k^{SU}$ of $SU_n$ with highest
  weight $(\lambda_i)_i$ such that $\sum\lambda_i = -l \in\{-n+1, \ldots,
  0\}$. Then the subobjects $w' \subset w\otimes v^{\otimes l}$ have highest
  weights $(\lambda'_i)_i$ such that $\sum\lambda'_i = 0$, hence they factor to
  representations of $PU_n$. Since we clearly have $w' \in B_{k+l}^{SU} \subset
  B_{k+n}^{SU}$, the preceeding discussion shows that $w' \in
  B_{a(k+n)}^{PU}$. But for any such $w'$ we also have $w \subset w' \otimes \bar
  v^{\otimes l}$. We have thus proved:
  $$B_k^{SU} \subset B_{a(k+n)}^{PU} \otimes U\quad\quad
  \text{where}\quad\quad U = \textstyle\bigoplus_{l=0}^{n-1} \bar v^{\otimes l}$$ 

  From this inclusion it clearly follows $b_k^{SU} \leq (\dim U)^2
  b_{a(k+n)}^{PU}$. Hence we have alltogether $b_k^{PU} \leq b_{2k}^A \leq
  b_{2k}^{SU} \leq (\dim U)^2 b_{a(2k+n)}^{PU}$. Since the $(b_k)$ sequences are
  growing and $(b_k^{SU})$ has growth exponent $d = n^2-1$, this proves that
  $(b_k^{PU})$, $(b_k^A)$ are polynomially growing with the same exponent.
\end{proof}

\begin{corollary}
  The discrete quantum group associated with $A_o^*(n)$ is amenable and has
  the Property of rapid decay.
\end{corollary}

\begin{proof}
  See \cite{bv1}, Proposition~2.1 and \cite{ve2}, Proposition~4.4.
\end{proof}

To illustrate the proof of Theorem~9.2 and the Cayley graph computations of
Section~8, let us draw the Cayley graphs of $C(PU_3)$, $A_o^*(3)$ and $C(SU_3)$
in $\bar X_{++}$. These graphs have no multiplicity and no loops, except in the
case of $PU_3$ where there are $2$ loops, that we do not represent, at each
vertex different from the origin. The arrows denote the root system of $SU(3)$,
and the dots are the images of $e_1$, $e_2$, $e_3$ in $\bar X_{++}$.

\begin{figure}[h!]
  \centering
  \includegraphics[scale=0.35]{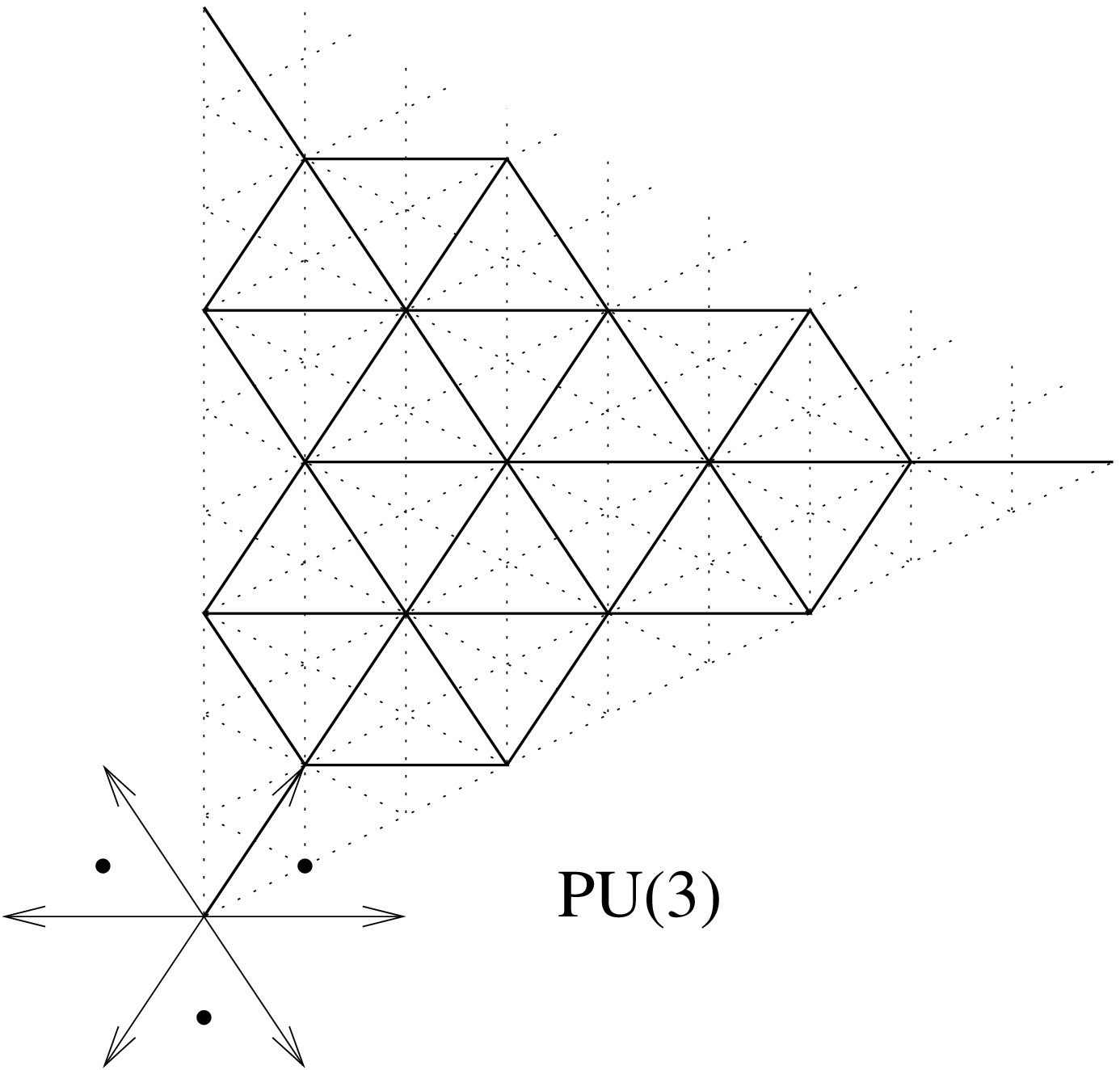}
  \includegraphics[scale=0.35]{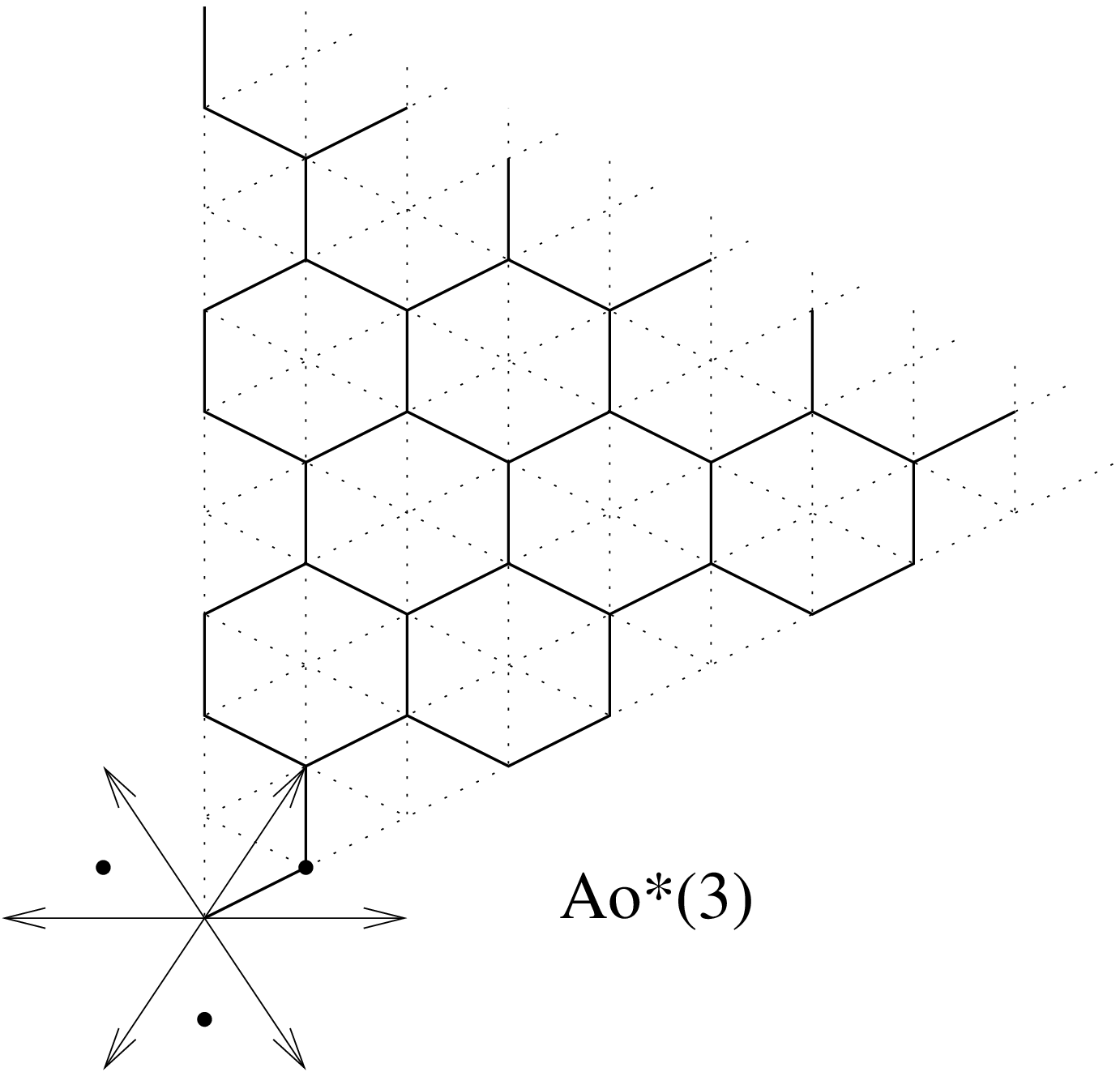}
  \includegraphics[scale=0.35]{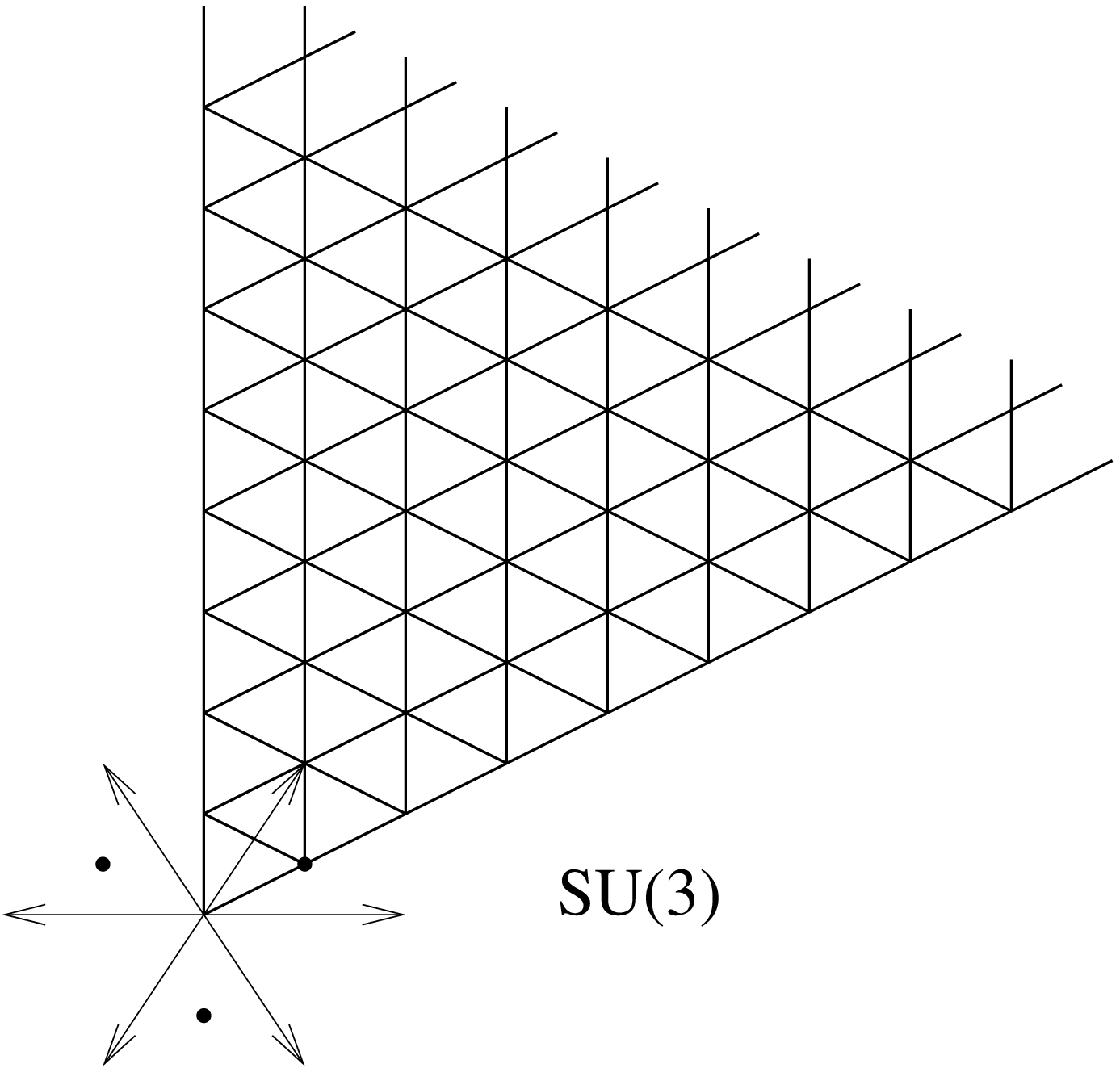}
\end{figure}

As a conclusion, let us remark that the results of this paper and the pictures
above seem to hint at the existence of some geometrical data behind compact
quantum groups, in the spirit of the classical constructions for semisimple Lie
algebras. It is tempting to ask whether the theoretical framework of \cite{wo4}
for differential calculus on compact quantum groups can give more insight on the
nature of the geometrical objects involved. However this seems to be a very
challenging question, and other intermediate examples between the classical
world and the free world of compact quantum groups might be needed first.

\end{document}